\documentclass[10pt,leqno]{amsart}
\usepackage{graphicx}
\usepackage{indentfirst,csquotes}

\usepackage{amssymb,amsthm,amsmath}
\usepackage{xcolor,paralist,hyperref,titlesec,fancyhdr,etoolbox}

\hypersetup{ colorlinks=true, linkcolor=black, filecolor=black, urlcolor=black }

\def\Card{\mathop{\rm Card}\nolimits}
\def\Char{\mathop{\rm Char}\nolimits}
\def\Munom{\mathop{\rm Munom}\nolimits}
\def\Dom{\mathop{\rm Dom}\nolimits}
\def\Fun{\mathop{\rm Fun}\nolimits}
\def\Op{\mathop{\rm Op}\nolimits}
\def\E{{\mathcal{E}}}
\def\N{{\mathbb{N}}}
\def\R{{\mathbb{R}}}

\begin{document}

\title{A Better Linear Unbiased Estimator for Averages over Discrete Structures}
\author[B J Braams]{Bastiaan J. Braams}
\date{\today}
\address{Centrum Wiskunde \& Informatica (CWI), Amsterdam, The Netherlands}
\email{b.j.braams@cwi.nl}

\begin{abstract}
Given an i.i.d.\ sample drawn from some probability distribution on a finite set, the best (in the sense of least variance) linear unbiased estimator (BLUE) of the average of any quantity with respect to that distribution is the sample average of the quantity.
Here we consider the situation in which, together with the sample, also the probability mass (possibly unnormalized) at each sample point is provided.
We show that with that information BLUE can be systematically improved.
The proposed procedure is expected to have applications in statistical physics, where it is common to have a closed-form specification of the relevant (unnormalized) probability distribution.
\end{abstract}

\maketitle

\let\thefootnote\relax

\bigskip

\section*{Introduction}

The present work is concerned with Monte Carlo evaluation of averages with respect to a probability distribution that is defined on a discrete space.
Let $D$ be that space and let $p:D\to\R_+$ be an unnormalized probability distribution with (unknown) normalization constant $Z$.
Using vector notation with vectors indexed over $D$, $Z=p^T.e$, where superscript $T$ denotes the transpose and $e$ is the all-ones vector.
We consider the situation, common in statistical physics, that a sampling procedure is available to generate $i\sim p/Z$ and that for each point in the sample the probability mass $p(i)$ is given or can be computed.
The result of a sampling process may be encoded in a count vector $c:D\to\N$, in which $c(i)$ is the number of times that element $i\in D$ is found in the sample.
We imagine that a hash table will be used to store the nonzero $c(i)$ together with the associated $p(i)$.

Let $f:D\to\R$ be a quantity of interest, and it is desired to estimate $\bar f=Z^{-1}p^T.f$.
The standard Monte Carlo estimate, here expressed using the count vector $c$, is $\bar f\simeq N^{-1}c^T.f$.
This estimator is linear in $f$ and is unbiased, because for each $i$, $c(i)/N$ is an unbiased estimate of $p(i)/Z$.
It is known by the acronym BLUE, the Best Linear Unbiased Estimator.

The question naturally arises, can we do better, and more specifically, can we do better within the class of linear unbiased estimators?
In a standard setting in applied statistics, only the counts $c(i)$ and function values $f(i)$ are available, but in our setting in statistical physics, also the $p(i)$ are known.
We will show that indeed, using all the data, a better linear unbiased estimator is available.

Given a count vector $c$ that is the result of sampling independently $N$ times from the distribution $p/Z$, define $S\subseteq D$ by $i\in S\Leftrightarrow1\leq c(i)$.
$S$ is a random set, and (given $N$, but before sampling) for any $i\in D$ the probability of $i\in S$ is $q(i)=1-(1-p(i)/Z)^N$: one minus the probability of obtaining $c(i)=0$.
It follows that the function $w$, defined on $D$ by $w(i)=p(i)/q(i)$ for $i\in S$ and $w(i)=0$ for $i\in D\setminus S$, is an unbiased estimate of $p$.
For $\bar f$ this leads to the estimator $\bar f\simeq Z^{-1}w^T.f$, which is linear in $f$ and unbiased.
The inner product reduces to a sum over $S$ only, just like the BLUE estimator.
Note that the quantity $Z$ is explicitly used here, and we will describe how to obtain it within the present context.
For now, assume that $Z$ is known.

The primary objective of the present work is to show that the variance of the estimator $\bar f\simeq Z^{-1}w^T.f$ is at most the variance of BLUE, up to a small correction associated with uncertainty in the normalization constant $Z$.
The reduction of variance is significant if a significant fraction of the total mass of the distribution is concentrated on points that are expected to be found at least once in the sample.
In the limit that all the mass is concentrated on such points, the proposed procedure tends smoothly to zero-variance exact summation.
At the other extreme, if the probability distribution has a very small fraction of the total mass on points that may be expected to occur at least once in the sample, then the proposed procedure reduces to BLUE.

\section*{Notation}

In the remainder of this manuscript we will employ some notational conventions associated primarily with E.\ W.\ Dijkstra~\cite{DS}.
If $\Op$ is any associative and commutative operation then $\langle\Op i,j:\Dom(i,j):\Fun(i,j)\rangle$ is the result of applying $\Op$ to the values $\Fun(i,j)$ as the indices $i,j$ vary over the domain indicated by the boolean expression $\Dom(i,j)$. This applies, for example, to sum ($\Sigma$), product ($\Pi$), Min, Max, set union and intersection, and existential and universal quantification.
If $f$ is a vector, matrix, function, or other structure that involves indices or arguments, and if $\Fun(f)$ is a boolean expression that may be understood pointwise, then $[\Fun(f)]$ means $\Fun(f)$ at all points of the domain of $f$. (Example: $[0\leq f]$ to say that $f$ is everywhere nonnegative.)

We review some notation that has already been introduced and that will retain its meaning.
$D$ is the fixed finite domain of interest.
It is expected to be made up of combinatorial objects and to be exponentially large in some natural parameter.
Elements of $D$ will be called points or objects, and we try to always use $i$ to refer to an element of $D$.
$p:D\to\R_+$ is a probability distribution (possibly unnormalized) on $D$.
$Z=\langle\Sigma i:i\in D:p(i)\rangle$ is the normalization constant.
$f:D\to\R$ is a quantity of interest; i.e., we want to estimate averages $\langle\Sigma i:i\in D:p(i)f(i)\rangle/\langle\Sigma i:i\in D:p(i)\rangle$.
In the context of sampling, $N$ is the number of sample points including repetitions.
Given $N$, but before drawing a sample, $q:D\to[0,1]$ is defined by $q(i)=1-(1-p(i)/Z)^N$; this is the probability that the point $i$ is included in a sample obtained by $N$ draws from $p/Z$.
After the draw, $c:D\to\N$ is the count vector associated with the sample.
$S$ is the set of points in $D$ included in the sample.
$M=\Card(S)$, the cardinality of $S$.
(Therefore, $M=\langle\Sigma i:i\in S:1\rangle$ and $N=\langle\Sigma i:i\in S:c(i)\rangle$.)
$e$ is the base of the natural logarithm or it is the all-ones vector; I think that the two meanings will not be confused.
For any subset $D'\subseteq D$ we use $P(D')$ for the total mass on $D'$: $P(D')=\langle\Sigma i:i\in D':p(i)\rangle$.
$\Char$ is the characteristic function of a relation ($\Char({\rm false})=0$, $\Char({\rm true})=1$).
$\Card$ is used for the cardinality of a set.

\section*{Determining the normalization constant}

If $Z$ is known, then the proposed method to evaluate $\bar f$ employs the approximation
\begin{equation}
    \bar f \simeq Z^{-1}\,\langle\Sigma i:i\in S: (p(i)/(1-(1-p(i)/Z)^N))\cdot f(i)\rangle
    \label{eq:avf}
\end{equation}
It is natural to require that this approximation is an identity when applied to the constant function $f$: $[f=1]$.
That requirement leads to a nonlinear equation for $Z$:
\begin{equation}
    Z = \phi(Z)\;,\quad \phi(Z)=\langle\Sigma i:i\in S: p(i)/(1-(1-p(i)/Z)^N)\rangle
    \label{eq:itz}
\end{equation}
In this section, we investigate properties of the iteration $Z:=\phi(Z)$.
We will show that (away from boundary cases that will be identified) the iteration map is a contraction on the range $P(S)\leq Z$, with $\phi$ increasing and convex.
This implies that simple Picard iteration, but also secant iteration and Newton iteration, are all suitable to determine $Z$ in eq.~\ref{eq:itz}.
We begin by dismissing the boundary cases.

If $N=0$ (there are no data), then $S=\emptyset$ and the iteration produces $Z=0$ immediately.
All averages are evaluated as $0^{-1}\cdot0$, which seems a fine result for this truly anomalous case.

If $N=1$ then $S$ is a singleton and the iteration takes on the form $Z:=Z$; in other words, $Z$ is entirely undetermined beyond the range requirement $P(S)\leq Z$.
However, eq.~\ref{eq:avf} is correctly defined; it says $\bar f=f(i)$ where $i$ is the unique element of $S$.
(We will see more generally that in cases where $Z$ is poorly defined by eq.~\ref{eq:itz}, nevertheless eq.~\ref{eq:avf} is well-defined.)

If $M=1$ ($S$ is a singleton) and $2\leq N$ then $\phi(Z)=P(S)/(1-(1-P(S)/Z)^N)$.
In this case the fixpoint is $Z=P(S)$; all mass is assumed to be concentrated on $S$ and the approximation eq.~\ref{eq:avf} becomes $\bar f=f(i)$ where $i$ is the unique element of $S$.
Also in this case, one may quibble with $Z=P(S)$, but $\bar f=f(i)$ is the only reasonable assignment.

If $M=N$ and $2\leq N$ ($[c=1]$ on $S$) then we use that $Z/N<p(i)/(1-(1-p(i)/Z)^N)$ to find $Z<\phi(Z)$ always.
The iteration will diverge with $Z\to\infty$.
Also in this limit eq.~\ref{eq:avf} is correctly defined; we obtain $\bar f=N^{-1}\langle\Sigma i:i\in S:f(i)\rangle$,
which is just the classical BLUE result.

We proceed to the main case, $2\leq M<N$, and as before we consider the iteration $Z:=\phi(Z)$, eq.~\ref{eq:itz}, on the range $P(S)\leq Z$.
Note that $[0<p]$ on $S$: points that have zero mass will not be sampled.
At $Z=P(S)$, $Z<\phi(Z)$.
As $Z\to\infty$, $\phi(Z)/Z\to M/N$, and so eventually $\phi(Z)<Z$.
We want to establish that $\phi$ is increasing, $0<\phi'(Z)$ for $P(S)<Z$, and that it is convex, $0<\phi''(Z)$ on the same domain.
Together with the asymptotic result $\phi(Z)\to(M/N)Z$ this will establish that $\phi$ is a strict contraction mapping (away from the $M$, $N$ boundary cases already discussed).

The sign properties $0<\phi'(Z)$ and $0<\phi''(Z)$ will be established by showing that they hold for each separate term in the sum in eq.~\ref{eq:itz}.
Define the auxiliary function $\psi(t)=1/(1-(1-t/N)^N)$ on the domain $0<t\leq N$.
The individual terms in eq.~\ref{eq:itz} are $\phi_i(Z)=p(i)/(1-(1-p(i)/Z)^N)$, $=p(i)\psi(Np(i)/Z)$.
We need to establish $0<-t^2\psi'(t)$ for $0<t<N$ and $0<t^2(t^2\psi'(t))'$ on the same domain.
The first inequality is easy, $-t^2\psi'(t)=t^2(1-t/N)^{N-1}/(1-(1-t/N)^N)^2$, and with $0<t<N$ the strict inequality is clear.
(It is an equality when $t=N$.)

The second-derivative inequality requires more work.

In order to guide our effort we first see how to establish the same second-derivative inequality for the related function $\eta(t)=1/(1-e^{-t})$ that is obtained in the limit $N\to\infty$.
We obtain $-t^2\eta'(t)=t^2e^{-t}/(1-e^{-t})^2$ and (with a bit of care in choosing the most appropriate factorization), $t^2(t^2\eta'(t))'=t^3e^{-t}(1-e^{-t})(2+t)(e^{-t}-(1-t/2)/(1+t/2))$.
The first four factors are manifestly positive and the final factor is recognized as the error term for the $[1,1]$ Pad\'e approximation to $e^{-t}$, which we know to be of definite sign.
Let us remind ourselves how to prove positivity of that final factor.
At $t=0$, $e^{-t}-(1-t/2)/(1+t/2)=0$.
Differentiate to find $(d/dt)(e^{-t}-(1-t/2)/(1+t/2))=-e^{-t}+(1+t/2)^{-2}$.
Then use $1+t/2<e^{t/2}$ to establish $0<(d/dt)(e^{-t}-(1-t/2)/(1+t/2))$, and thereby establish $0<t^2(t^2\eta'(t))'$ for $0<t$.

Returning to the original function $\psi$, we start with $-t^2\psi'(t)=t^2(1-t/N)^{N-1}/(1-(1-t/N)^N)^2$, $=t^2(1-t/N)^{N-1}(\psi(t))^2$.
Next, with guidance from the previous effort on the second-derivative inequality for $\eta$ and with judicious choice of factorization, we obtain:
\begin{equation}
\begin{split}
    t^2(t^2\psi'(t))' &= t^3\cdot(-\psi'(t))\cdot\psi(t)\cdot(1-t/N)^{-1}\cdot(2+t-t/N)\\
    &\;\cdot\left((1-t/N)^N-\frac{2-t-t/N}{2+t-t/N}\right)
\end{split}
\label{eq:psipp}
\end{equation}
All factors except the final one are obviously positive.
The boundary cases $N=1$ and $t=N$ provide a sanity check on eq.~\ref{eq:psipp}. (The singularity of $(1-t/N)^{-1}$ at $t=N$ is cancelled in $\psi'(t)$.)
For the final factor, as in dealing with the $N\to\infty$ limit function $\eta$, we observe that at $t=0$, $(1-t/N)^N-(2-t-t/N)/(2+t-t/N)=0$, and we differentiate the expression to obtain $-(1-t/N)^{N-1}+(1+(1-N^{-1})t/2)^{-2}$.
To show that this is positive we observe $1+(1-N^{-1})t/2<(1-t/N)^{(1-N)/2}$.

This concludes the demonstration that in the original iteration, eq.~\ref{eq:itz}, the right hand side function $\phi$ is increasing and (strictly) convex away from the identified boundary cases $M=1$ and $M=N$.
The weights in eq.~\ref{eq:avf} are always well defined and the entire procedure has correct limits.
If all the mass is concentrated on $S$ and $N\to\infty$, then the limit is the plain weighted sum $\bar f=(P(S))^{-1}\langle\Sigma i:i\in S:p(i)f(i)\rangle$.
If $M=N$ then the traditional BLUE estimator is recovered.

As a starting guess for the iteration of eq.~\ref{eq:itz} the value $Z=P(S)$ can be used, but a much better starting guess is provided by the Good-Turing~\cite{GT} missing mass estimator.
Let $M'=\langle\Sigma i:c(i)=1:1\rangle$; so $M'\leq M\leq N$.
The Good-Turing missing mass estimator is $1-P(S)/Z=M'/N$.
It is undefined ($Z\to\infty$) if $M'=N$, but it does not matter; this is the known special case $M=N$.

\section*{Calculation of means and covariances}

The mean and the covariance of the BLUE estimator are of course very well known, but we repeat the calculation here in order to have everything together for a comparison with the mean and covariance of the estimator proposed in this manuscript.

To assess BLUE we require the expectation values $\E(c(i)/N)$ and $\E(c(i)c(j)/N^2)$ where the count vector $c$ is a random variable obtained by sampling $N$ times from the distribution $p/Z$.
We employ multinomial coefficients with the notation $\Munom(c)=(c^T.e)!/c!$, $=N!/c!$, where $c!$ is shorthand for $\langle\Pi i:i\in D:c(i)!\rangle$ (and $N=c^T.e$ as before).
The fundamental identity for multinomial coefficients is $\langle\Sigma c:c^T.e=N:\Munom(c)p^c\rangle=(p^T.e)^N$, in which $p^c$ is shorthand for $\langle\Pi i:i\in D:p(i)^{c(i)}\rangle$, with $0^0=1$.
(All sums and products over $D$ could be restricted to $S$ without changing the value.)
To review:
\begin{equation}
\begin{split}
    \E(c(i)/N)&=\langle\Sigma c:c^T.e=N:\Munom(c)(p/Z)^c(c(i)/N)\rangle\\
    &=N^{-1}Z^{-N}p(i)(\partial/\partial p(i))\langle\Sigma c:c^T.e=N:\Munom(c)p^c\rangle\\
    &=N^{-1}Z^{-N}p(i)(\partial/\partial p(i))(p^T.e)^N\\
    &=Z^{-N}p(i)(p^T.e)^{N-1}\;,\quad=p(i)/Z
\end{split}
\end{equation}
And for the covariance matrix elements:
\begin{equation}
\begin{split}
    &\E(c(i)c(j)/N^2)-p(i)p(j)/Z^2\\
    &\qquad=p(i)(\partial/\partial p(i))p(j)(\partial/\partial p(j))(p^T.e/Z)^N-p(i)p(j)/Z^2\\
    &\qquad=(1-N^{-1})p(i)p(j)/Z^2+\delta(i,j)N^{-1}p(i)/Z-p(i)p(j)/Z^2\\
    &\qquad=N^{-1}\left(\delta(i,j)(p(i)/Z)(1-p(i)/Z)-(1-\delta(i,j))p(i)p(j)/Z^2\right)
\end{split}
\label{eq:covblue}
\end{equation}
Sanity check: $\langle\Sigma i,j::\E(c(i)c(j)/N^2)-p(i)p(j)/Z^2\rangle=0$; BLUE is exact on the constant function $[f=1]$.

Now the related calculation for the weights proposed in this work and assuming that $Z$ is known.

Expectation values:
\begin{equation}
\begin{split}
    \E(w(i)/Z)&=Z^{-1}\left\langle\Sigma c:c^T.e=N:\Munom(c)(p/Z)^c\Char(1\leq c(i))p(i)/q(i)\right\rangle\\
    &=Z^{-1}(p(i)/q(i))(1-(1-p(i)/Z)^N),\quad=p(i)/Z
\end{split}
\end{equation}
(using $\Char(1\leq c(i))=1-\Char(c(i)=0)$).

For the covariance matrix elements a similar computation:
\begin{equation*}
\begin{split}
    &\E(w(i)w(j))\\
    &\qquad=\left\langle\Sigma c:c^T.e=N:\Munom(c)p^c\Char(1\leq c(i)\wedge1\leq c(j))p(i)p(j)/(q(i)q(j))\right\rangle\\
    &\qquad=(p(i)p(j)/(q(i)q(j)))\\
    &\qquad\quad\cdot\left(1-(1-p(i)/Z)^N-(1-p(j)/Z)^N\right.\\
    &\qquad\qquad\quad\left.+\,\delta(i,j)(1-p(i)/Z)^N+(1-\delta(i,j))(1-p(i)/Z-p(j)/Z)^N)\right)
\end{split}
\end{equation*}
(using $\Char(1\leq c(i)\wedge1\leq c(j))=1-\Char(c(i)=0)-\Char(c(j)=0)+\Char(c(i)=0\wedge c(j)=0)$, and splitting up the final $\Char()$ into cases $i=j$ and $i\neq j$).
Then:
\begin{equation}
\begin{split}
    &\E(w(i)w(j)/Z^2)-p(i)p(j)/Z^2\\
    &\qquad=\delta(i,j)(p(i)/Z)^2(1/q(i)-1)\\
    &\qquad-(1-\delta(i,j))((p(i)p(j)/Z^2)/(q(i)q(j)))\\
    &\qquad\qquad\cdot\left((1-p(i)/Z-p(j)/Z+p(i)p(j)/Z^2)^N-(1-p(i)/Z-p(j)/Z)^N\right)
\end{split}
\label{eq:covnew}
\end{equation}

Now I would like to assert that the covariance matrix in eq.~\ref{eq:covnew} is dominated (in the positive semidefinite matrix ordering) by that of eq.~\ref{eq:covblue}, but this will not fly.
The BLUE procedure is exact for constant vectors.
In the new procedure proposed here, if the normalization constant is taken to be known before the sampling, then for constant vectors we obtain the correct expectation value, but with a nonzero uncertainty.
And if the normalization constant is determined after sampling then eq.~\ref{eq:covnew} is not strictly valid.
I expect this to be a very small effect in practice, but it certainly confuses the analysis.

So let us focus attention on the diagonal elements only.
From eq.~\ref{eq:covblue}, $\E(c(i)^2/N^2)-p(i)^2/Z^2=(1/N)(p(i)/Z)(1-p(i)/Z)$.
From eq.~\ref{eq:covnew}, $\E(w(i)^2/Z^2)-p(i)^2/Z^2=(p(i)/Z)^2(1/q(i)-1),\,=(p(i)/Z)^2(1-p(i)/Z)^N/(1-(1-p(i)/Z)^N)$.
If $Np(i)\ll1$ (point $i$ is unlikely to be in the sample) then $\E(w(i)^2/Z^2)-p(i)^2/Z^2\simeq(p(i)/Z)^2(1-Np(i)/Z)/(Np(i)/Z),\;=(1/N)(p(i)/Z)(1-Np(i)/Z)$, which is essentially the BLUE result.
On the other hand, if $1\ll Np(i)/Z$ then $(1-p(i)/Z)^N\simeq\exp(-Np(i)/Z)$, so the contribution of those points to the total variance vanishes exponentially.

In conclusion, as quantified above, the proposed procedure is expected to be valuable if a significant fraction of the total probability mass is concentrated on points that are expected to be sampled more than once.

\section*{Numerical tests}

Fortran codes for numerical tests of the estimator of this manuscript are being assembled in the GitHub location {\tt https://github.com/bjbraams/sampling}.
The tests are synthetic and for this purpose the particular combinatorial structure of the space $D$ is irrelevant.
Only the density of states matters.
For the first tests, $D$ contains 96-bit unsigned integers encoded into three default integers.
Each $i\in D$ represents the real number $x=i+1/2$ and we create a distribution $p(i)$ that decays as $(b/a)(1+x/a)^{-b-1}$, but for precision it is obtained by finite differences of an associated cumulative distribution.
All arithmetic on the values $x$ represented by $i\in D$ is done via 64-bit floating point arithmetic on the three integer values that encode $i$.
Parameters $a$, $b$ and the sample size $N$ are set to explore different regimes: most of the mass expected to be sampled, most of the mass in the tail, or anything in between.
For the test function $f$ some highly oscillatory function is used that (as a continuous function over the associated continuous density) would integrate to zero.

Tests run so-far are consistent with the expected behaviour of the sampling method.
This will be described further.

\section*{Closing remarks}

At first sight the proposed procedure requires saving all unique states $i$ in the sample.
But it suffices to save only some cryptographic hash of the state (to avoid double counting) together with the quantity of interest.
And if even the quantity of interest is too large to save every instance $f(i)$, it can be accumulated with weights belonging to some values of $Z$ and $N$, and finally interpolated when $Z$ and $N$ are known.

In order to be effective the proposed procedure requires some concentration of measure, and as presented here it would be completely useless for continuous domain $D$.
This has to be addressed by some kind of binning, assigning points to cells, but then one needs estimates for the mass of each cell.
This seems to require special analysis for any particular class of problems.
(I have tried to find a general procedure involving only the density of states, but without success.)

\section*{Acknowledgments}

I thank Wouter Koolen and Peter Gr\"unwald of CWI for valuable conversations.

\end{document}